\documentclass[11pt]{article}

\newcommand{\documentdate}{July 25th, 2017}

\usepackage{a4wide,latexsym,mystyle,varioref,amsmath}
\newcommand{\eqn}[2]{\begin{equation}\label{#1}{#2}\end{equation}}
\newcommand{\disp}[1]{\[{#1}\]}

\newcommand{\ieqo}{i=1,2}
\newcommand{\lambdaleft}{\lambda_{\hbox{\scriptsize left}}}
\newcommand{\flow}{f_{\rm low}}

\title{Improved second-order evaluation complexity for unconstrained
  nonlinear optimization using high-order regularized models}

\author{
C. Cartis\thanks{Mathematical Institute,
   Oxford University,
   Oxford OX2 6GG, England.  Email: coralia.cartis@maths.ox.ac.uk},
N. I. M. Gould\thanks{ Computational Mathematics Group,
   STFC-Rutherford Appleton Laboratory,
   Chilton OX11 0QX, England. Email:  nick.gould@stfc.ac.uk .
 The work of this author was supported by EPSRC grant EP/M025179/1}
~and~Ph. L. Toint\thanks{ Namur Center for Complex Systems (naXys) and
  Department of Mathematics, University of Namur, 61, rue de
  Bruxelles, B-5000 Namur, Belgium.  Email: philippe.toint@unamur.be}
}

\date{\documentdate}

\begin{document}

\maketitle

\begin{abstract}
The unconstrained minimization of a sufficiently smooth objective
function $f(x)$ is considered, for which derivatives up to order $p$, $p\geq 2$, are assumed to be available. An adaptive regularization algorithm is proposed 
that uses Taylor models of the objective of order $p$ and that  is guaranteed
to find a first- and second-order critical point in at most
$O \left(\max\left( \epsilon_1^{-\frac{p+1}{p}}, \epsilon_2^{-\frac{p+1}{p-1}}
\right) \right)$
function and derivatives evaluations, where $\epsilon_1$ and $\epsilon_2 >0$
are prescribed first- and second-order optimality tolerances. 
Our approach extends the method in Birgin et al. (2016) to finding second-order critical points, and establishes the novel complexity bound for second-order criticality
under identical problem assumptions as for first-order, namely, that the $p$-th derivative tensor is
Lipschitz continuous and that $f(x)$ is bounded from below.
The evaluation-complexity
bound for second-order criticality improves on all such  known existing results. 
\end{abstract}

\numsection{Introduction}
A question of general interest in computational optimization is to know
how many evaluations of the functions that define a given problem are needed
for an algorithm to find an estimate of a local minimizer.
Considerable advances have been made on this topic, both for convex problems
\cite{Nest04} and nonconvex ones \cite{CartGoulToin12:opt}.
Although much of this research has been devoted to the important issue of
finding approximate first-order critical points, some authors have addressed
the case where higher-order necessary optimality conditions must also be
satisfied.

We consider the unconstrained minimization of a $C^2$
objective function $f:\Re^n \rightarrow \Re$. It is, of course,
well known that a finite minimizer $x_*$ of $f$ necessarily satisfies
the first- and second-order criticality conditions $\nabla_x f(x_*) = 0$ and
$\lambdaleft(\nabla^2 f_2(x_*) )\geq 0$,
where $\lambdaleft$ denotes the leftmost eigenvalue of its symmetric
matrix argument. Thus a reasonable requirement might be to find a point $x_k$
for which
\eqn{opt12}{\|\nabla_x f(x_k)\| \leq \epsilon_1 \tim{and}
\lambdaleft\left(\nabla^2_x f(x_k)\right) \geq - \epsilon_2}
for given, small $\epsilon_1, \epsilon_2 > 0$ and suitable norm $\|\cdot\|$.

The earliest analysis we are aware of that provides both first- and second-order evaluation complexity guarantees
considers  cubic
regularization  methods and shows that at most
\eqn{arc-eval}{O\left(\max\left(\epsilon_1^{-3/2},\epsilon_2^{-3}\right)\right)}
evaluations of $f$ are required to satisfy \req{opt12}
so long as the objective function is bounded from below,
and its Hessian is Lipschitz continuous \cite{NestPoly06}. Adaptive cubic regularization variants with inexact subproblem solves
and similar guarantees were proposed in 
\cite{CartGoulToin11b:mp, CartGoulToin12d}.
Under similar conditions, many trust-region (TR) algorithms require at most
$O\left(\max\left(\epsilon_1^{-2},\epsilon_2^{-3}\right)\right)$
evaluations. Crucially, examples are
known for which such order estimates are tight both for trust-region and regularization methods \cite{CartGoulToin12d}.
Of late, more sophisticated trust region methods and quadratic regularization ones have been
proposed that echo the order of the ARC estimates
\cite{CurtRobiSame17:mp,MartRayd16, BirginMartinez17}. At the same time, other methods
\cite{GrapYuanYuan16:jota,GratRoyoVice17} have been shown to mirror
the TR-like evaluation estimate in a more general or simplified way, respectively.

The fact that the best-known evaluation bound for ARC is essentially tight, suggests
that in order to do better, one needs to add further ingredients. A
similar picture emerged for evaluation bounds for first-order critical
points: improved bounds of order $\mathcal{O}\left(\epsilon_1^{-\frac{p+1}{p}}\right)$, $p\geq 2$,  were obtained in \cite{BirgGardMartSantToin17:mp}
for $p$-times continuously differentiable functions using regularization methods that employ
higher-order local models.
This will be the theme
here. In order to improve upon the estimate \req{arc-eval} for second-order criticality, we will use
a higher-order model and regularization. The model minimization
conditions however, are approximate and local, for both first- and
second-order criticality.

In \S\ref{sec-algorithm}, we define terminology and propose our
new algorithm, while in \S\ref{sec-analysis}, we provide a convergence
analysis that indicates an improved complexity bound. We provide
further comments and perspectives in \S\ref{sec-conclusion}.


\numsection{A regularized $p$-th order model and algorithm}
\label{sec-algorithm}

Let $p\geq 2$. Consider the optimization problem
\beqn{problem}
\min_{x \in \smallRe^n}  f(x),
\eeqn
where we assume that 
$f\in  \mathcal{C}^{p,1}(\Re^n)$, namely, that:
\begin{itemize}
\item  $f$ is
$p$-times continuously differentiable;

\item $f$ is bounded below by $f_{\rm low}$

\item  the $p$-th derivative of $f$ at $x$, the $p$-th order tensor
\[
\nabla_x^pf(x) = \left[ \frac{\partial^p f(x)}{\partial x_{i_1}\ldots\partial
    x_{i_p}} \right]_{i_j \in \ii{n}, j=1, \ldots,p}
\]
is globally Lipschitz continuous, that is, 
there exists a constant $L \geq 0$ such
that, for all $x,y \in \Re^n$,
\beqn{tensor-Lip}
\| \nabla_x^p f(x) - \nabla_x^p f(y) \|_{[p]} \leq (p-1)!\, L \| x-y\|.
\eeqn
\end{itemize}
In \req{tensor-Lip}, $\|\cdot\|_{[p]}$ is the tensor norm recursively induced
by the Euclidean norm $\|\cdot\|$ on the space of $p$-th order tensors, which
is given by
\beqn{Tnorm}
\|T\|_{[p]} \eqdef \max_{\|v_1\|=\cdots=\|v_p\|=1}  | T[v_1,\ldots,v_p] |,
\eeqn
where $T[v_1,\ldots,v_j]$ stands for the tensor of order $p-j \geq 0$ resulting from the
application of the $p$-th order tensor $T$ to the vectors
$v_1,\ldots,v_j$\footnote{Note that $\|\cdot\|_{[1]}=\|\cdot\|$, the usual Euclidean vector norm.}. Let $T_p(x,s)$ be the
Taylor series of the function $f(x+s)$ at $x$ truncated at order~$p$
\beqn{Taylor-def}
T_p(x,s) \eqdef f(x) + \sum_{j=1}^p \frac{1}{j!}\nabla_x^jf(x)[s]^j,
\eeqn
where the notation $T[s]^j$ stands for the tensor $T$ applied $j$ times to the
vector $s$.

We shall use the following crucial bounds.

\llem{taylor-bounds-lemma}{[See Appendix~\ref{A-Lemma-proof}].
Let $f \in C^{p,1}(\Re^n)$, and $T_p(x,s)$ be the
Taylor approximation of $f(x+s)$ about $x$.
Then for all $x,s \in \Re^n$,
\eqn{resf}{f(x+s) \leq T_p(x,s) + \bigfrac{L}{p} \, \|s\|^{p+1},}
\eqn{resg}{\| \nabla_x^1 f(x+ s) -  \nabla_s^1 T_p(x,s) \|_{[1]} \leq L \|s\|^p}
and
\eqn{H-Lip}{
\| \nabla_x^2 f(x+ s) -  \nabla_s^2 T_p(x,s) \|_{[2]} \leq (p-1) L \|s\|^{p-1}.}
}

In order to describe our algorithm, we  define the regularized
Taylor series model
\beqn{mdef}
m(x,s, \sigma) = T_p(x,s) + \frac{\sigma}{p+1}\|s\|^{p+1},
\eeqn
whose gradient and Hessian are
\beqn{gradm}
\nabla_s^1 m(x,s,\sigma) = \nabla_s^1 T_p(x,s) +
\sigma \|s\|^p \frac{s}{\|s\|}
\eeqn
and
\beqn{Hm}
\nabla_s^2 m(x,s,\sigma) = \nabla_s^2 T_p(x,s) +
 \frac{\sigma}{p+1} \nabla_s^2 \left(\|s\|^{p+1}\right),
\eeqn
where
\beqn{Hr}
\nabla_s^2 \left(\|s\|^{p+1} \right)= (p+1)
\left[ (p-1) \|s\|^{p-3} s s^T + \|s\|^{p-1} I \right].
\eeqn
Note that
\beqn{atxk}
m(x,0,\sigma) = T_p(x,0) = f(x).
\eeqn

For the objective function $f$, we define
first- and second-order
criticality measures as
\beqn{chi1-def}
\chi_{f,1}(x) \eqdef \|\nabla_x^1f(x)\|
\eeqn
and
\beqn{chi2-def}
\chi_{f,2}(x)
\eqdef \max\Big[0,-\lambda_f(x)\Big] = \max\Big[0,- \min_{\|y\|=1}\nabla_x^2 f(x)[y]^2\Big]
\eeqn
where $\lambda_f(x) \eqdef \lambdaleft[\nabla_x^2 f(x)]$.
Similarly, for the model \req{mdef}, we consider the measures
\beqn{chim1-def}
\chi_{m,1}(x,s,\sigma) \eqdef \|\nabla_s^1m(x,s,\sigma)\|
\eeqn
and
\beqn{chim2-def}
\chi_{m,2}(x,s,\sigma)
\eqdef \max\Big[0,-\lambda_m(x,s,\sigma)\Big]
=\max\Big[0, - \min_{\|y\|=1}\nabla_s^2 m(x,s,\sigma)[y]^2\Big]
\eeqn
where $\lambda_m(x,s,\sigma) \eqdef \lambdaleft[\nabla_s^2 m(x,s,\sigma)]$.

\vspace*{2mm}
The minimization algorithm we consider is  now described in detail in
Algorithm 2.1 below. Note that if the second-order conditions are removed -- namely, the conditions for $i=2$ in \req{f-term} and \req{mterm} -- then this method reduces to the ARp algorithm in \cite{BirgGardMartSantToin17:mp}.\\

{\bf Algorithm 2.1: AR$p$.} 
\begin{description}
\item[Step 0: Initialization.]
  An initial point $x_0$  and an initial regularization parameter $\sigma_0>0$
  are given, as well as an accuracy levels $\epsilon_1$, $\epsilon_2$ and
  $\epsilon_3$.  The constants $\theta$,
  $\eta_1$, $\eta_2$, $\gamma_1$, $\gamma_2$, $\gamma_3$ and $\sigma_{\min}$
  are also given and satisfy
\beqn{eta-gamma2}
\theta > 0,  \ms
    \sigma_{\min} \in (0, \sigma_0], \ms
0 < \eta_1 \leq \eta_2 < 1 \tim{and}
0< \gamma_1 < 1 < \gamma_2 < \gamma_3.
\eeqn
Compute $f(x_0)$ and set $k=0$.

\item[Step 1: Test for termination. ]
Evaluate $\{\nabla_x^i f(x_k)\}_{i=1}^2$.
If
\beqn{f-term}
\chi_{f,i}(x_k) \leq \epsilon_i \tim{ for } \ieqo,
\eeqn
terminate with the approximate solution $x_\epsilon=x_k$. Otherwise compute
derivatives of $f$ from order 3 to $p$ at $x_k$.

\item[Step 2: Step calculation. ]
Compute the step $s_k$ by  approximately minimizing the model
$m(x_k,s,\sigma_k)$ with respect to $s$ in the sense that the conditions
\beqn{descent2}
m(x_k,s_k,\sigma_k) < m(x_k,0,\sigma_k)
\eeqn
and
\beqn{mterm}
\chi_{m,i}(x_k,s_k,\sigma_k) \leq \theta \|s_k\|^{p+1-i}, \ms (\ieqo)
\eeqn
hold.

\item[Step 3: Acceptance of the trial point. ]
Compute $f(x_k+s_k)$ and define
\beqn{rhokdef2}
\rho_k = \frac{f(x_k) - f(x_k+s_k)}{T_p(x_k,0) - T_p(x_k,s_k)}.
\eeqn
If $\rho_k \geq \eta_1$, then define
$x_{k+1} = x_k + s_k$; otherwise define $x_{k+1} = x_k$.

\item[Step 4: Regularization parameter update. ]
Set
\beqn{sigupdate2}
\sigma_{k+1} \in \left\{ \begin{array}{ll}
{}[\max(\sigma_{\min}, \gamma_1\sigma_k), \sigma_k ]  & \tim{if} \rho_k \geq \eta_2, \\
{}[\sigma_k, \gamma_2 \sigma_k ]          &\tim{if} \rho_k \in [\eta_1,\eta_2),\\
{}[\gamma_2 \sigma_k, \gamma_3 \sigma_k ] & \tim{if} \rho_k < \eta_1.
  \end{array} \right.
\eeqn
Increment $k$ by one and go to Step~1 if $\rho_k\geq \eta_1$ or to Step~2 otherwise.
\end{description}

Each iteration of this algorithm requires the approximate minimization of
$m(x_k,s,\sigma_k)$, and we note that conditions \req{descent2} and
\req{mterm} are always achievable as they are satisfied
at a
second-order
critical point of $m(x,s,\sigma)$.  Indeed, existing algorithms, such as the
standard second-order trust-region method \cite[\S6.6]{ConnGoulToin00}
and ARC \cite{CartGoulToin11:mp}
will find such a point as the regularized Taylor model is both sufficiently
smooth and bounded from below.\footnote{When $p$ is even, $m(x,s,\sigma)$
is smooth everywhere but at the origin, but a step from $s=0$
in the steepest-descent/eigen direction will move to a region for which
the model is always smooth.}
Moreover,
\emph{this approximate minimization does not involve additional
computations of $f$ nor its derivatives at points other than $x_k$, and
therefore the precise method used, and the resulting effort spent, in Step~2 have
no impact on the evaluation complexity}\footnote{We implicitly assume here that derivatives at $x_k$ can be stored explicitly.}.  Finally note that the second condition in
\req{mterm} disappears if $\lambdaleft(\nabla_x^2 T_p(x,s))\geq 0$.

Iterations for which $\rho_k \geq \eta_1$ (and hence $x_{k+1}=x_k+s_k$) are
called ``successful'' and we denote by $\calS_k \eqdef \{ 0 \leq j \leq k \mid
\rho_j \geq \eta_1 \}$  the index set of all successful iterations between 0
and $k$.  
We also denote the complement, $\calU_k$, of $\calS_k$ in
$\{0, \ldots, k\}$, that corresponds to the index set of ``unsuccessful''
iterations between 0 and $k$.  Note that, before termination, each successful
iteration requires the evaluation of $f$ and its first $p$ derivatives,
while only the evaluation of $f$ is needed at unsuccessful ones.

\numsection{Complexity analysis}
\label{sec-analysis}

As it is typical  for a complexity analysis of (regularization and other) methods, we proceed by showing lower bounds on the Taylor model decrease and on the length of the step at each iteration. 
The proofs of the next three lemmas is very similar to corresponding results in \cite{BirgGardMartSantToin17:mp} and hence we defer the proofs to the appendix (but still include them
for completeness, as the algorithm has changed).

\llem{Dm_lemma}{
The mechanism of Algorithm~2.1 guarantees that, for all $k  \geq 0$,
\beqn{Dphi}
T_p(x_k,0) - T_p(x_k,s_k) \geq \frac{\sigma_k}{p+1} \|s_k\|^{p+1},
\eeqn
and so \req{rhokdef2} is well-defined.
}

We next deduce a simple upper bound on the regularization parameter $\sigma_k$.

\llem{sigmaupper_lemma}{
Let $f\in \mathcal{C}^{p,1}(\Re^n)$.
Then, for all $k\geq 0$,
\beqn{sigmaupper}
\sigma_k
\leq \sigma_{\max} \eqdef \max\left[ \sigma_0,\frac{\gamma_3 L(p+1)}{p\,(1-\eta_2)}\right].
\eeqn
}

Our next move, very much in the line of the theory proposed in
\cite{CartGoulToin11b:mp, BirgGardMartSantToin17:mp}, is to show that the step cannot be
arbitrarily small compared with the gradient of the objective function
at the trial point $x_k+s_k$.

\llem{longs-g-lemma}{
Let $f\in \mathcal{C}^{p,1}(\Re^n)$.
Then, for all  $k\geq 0$,
\beqn{longs-g}
\|s_k\| \geq \left(\frac{\chi_{f,1}(x_k+s_k)}{L + \theta +\sigma_k}\right)^{\frac{1}{p}}.
\eeqn
}

Next we show that the step cannot also be arbitrarily small compared to the second order criticality measure \req{chi2-def} at the trial point $x_k+s_k$.
This is the crucial novel  ingredient of the paper, that is essential to the improved second-order complexity results.

\llem{longs-H-lemma}{
Let $f\in \mathcal{C}^{p,1}(\Re^n)$. Then, for
all  $k\geq 0$,
\beqn{longs-H}
\|s_k\| \geq \left(\frac{\chi_{f,2}(x_k+s_k)}{(p-1)L + \theta +
 p \sigma_k}\right)^{\frac{1}{p-1}}.
\eeqn
}

\proof{
Using \req{mdef} and the fact that $\min_z[a(z)+b(z)] \geq \min_z[a(z)]
+ \min_z[b(z)]$,
we find that
\disp{\arr{rl}{
\lambda_f(x_k+s_k) & \!\!\!\! =
\bigmin_{\|y\|=1}\nabla_x^2 f(x_k+s_k)[y]^2 \\
& \!\!\!\! = \bigmin_{\|y\|=1}
 \left( \nabla_x^2 f(x_k+s_k) - \nabla_s^2 T_p(x_k,s_k) - \frac{\sigma_k}{p+1}
\nabla_s^2 \|s_k\|^{p+1} + \nabla_s^2 m(x_k,s_k,\sigma_k) \right) [y]^2
 \\
& \!\!\!\! \geq \bigmin_{\|y\|=1}
 \left( \nabla_x^2 f(x_k+s_k) - \nabla_s^2 T_p(x_k,s_k) \right) [y]^2 +
\bigfrac{\sigma_k}{p+1}
\bigmin_{\|y\|=1} \left( - \nabla_s^2 \|s_k\|^{p+1} \right) [y]^2 + \\
& \bigmin_{\|y\|=1}
\nabla_s^2 m(x_k,s_k,\sigma_k) [y]^2
}}
Considering each term in turn, and using
\req{Tnorm} and \req{H-Lip}, we see that
\disp{\arr{rl}{
\bigmin_{\|y\|=1} & \!\!\!\!
 \left( \nabla_x^2 f(x_k+s_k) - \nabla_s^2 T_p(x_k,s_k) \right) [y]^2 \\
& \geq
\bigmin_{\|y_1\|=\|y_2\|=1}
 \left( \nabla_x^2 f(x_k+s_k) - \nabla_s^2 T_p(x_k,s_k) \right) [y_1,y_2] \\
& \geq -
\bigmax_{\|y_1\|=\|y_2\|=1} \left | \left( \nabla_x^2 f(x_k+s_k)
 - \nabla_s^2 T_p(x_k,s_k) \right) [y_1,y_2] \right | \\
& = - \| \nabla_x^2 f(x_k+s_k)  - \nabla_s^2 T_p(x_k,s_k) \|_{[2]}^{} \\
& \geq - (p-1)L\|s_k\|^{p-1},
\\*[2.5ex]
}}
and using \req{Hr}, we find that $\nabla_s^2 \left(\|s_k\|^{p+1} \right) [y]^2= (p+1)[(p-1)\|s_k\|^{p-3}(s_k^Ty)^2+\|s_k\|^{p-1}\|y\|^2]$, and so
\disp{
\bigmin_{\|y\|=1}   \left( - \nabla_s^2 (\|s_k\|^{p+1}) \right) [y]^2
 = - \bigmax_{\|y\|=1} \nabla_s^2 (\|s_k\|^{p+1}) [y]^2 =
 - p (p+1) \|s_k\|^{p-1}.
}
Recalling \req{chim2-def}, we have $\min_{\|y\|=1} \nabla_s^2 m(x_k,s_k,\sigma_k) [y]^2 =
\lambda_{m}(x_k,s_k,\sigma_k)$. This, and 
the last two displayed equations imply that
\eqn{interm}{-\lambda_f(x_k+s_k) \leq (p-1)L\|s_k\|^{p-1} + p \sigma_k \|s_k\|^{p-1} -\min[0, \lambda_{m}(x_k,s_k,\sigma_k)].}
As the right hand side of \req{interm} is nonnegative, the bound \req{interm} can be re-written as
\disp{
\max[0, -\lambda_f(x_k+s_k)]\leq  \left[(p-1)L + p \sigma_k\right] \|s_k\|^{p-1} +\max[0, -\lambda_{m}(x_k,s_k,\sigma_k)].
}
Combining the above with \req{chi2-def} and \req{chim2-def},  and with \req{mterm} with $i=2$, we conclude
\disp{
\arr{rl}{\chi_{f,2}(x_k+s_k) & \!\!\!\!\leq 
( (p-1)L + p \sigma_k )  \|s_k\|^{p-1} + \chi_{m,2}(x_k,s_k,\sigma_k)
\\ & \!\!\!\! \leq
( (p-1)L + \theta + p \sigma_k )  \|s_k\|^{p-1}}}
and \req{longs-H} follows.
}



We now bound the number of unsuccessful iterations as a function of the number
of successful ones and include a proof in the Appendix.

\llem{SvsU}{
\cite[Theorem~2.1]{CartGoulToin11b:mp}
The mechanism of Algorithm~2.1 guarantees that, if
\beqn{sigmax}
\sigma_{k} \leq \sigma_{\max},
\eeqn
for some $\sigma_{\max} > 0$, then
\beqn{unsucc-neg}
k +1 \leq |\calS_k| \left(1+\frac{|\log\gamma_1|}{\log\gamma_2}\right)+
\frac{1}{\log\gamma_2}\log\left(\frac{\sigma_{\max}}{\sigma_0}\right).
\eeqn
}

Using all the above results, we are now in position to state our main evaluation complexity result.

\lthm{final_theorem}{
Let $f\in \mathcal{C}^{p,1}(\Re^n)$.
Then, given $\epsilon_1> 0$ and $\epsilon_2>0$,
Algorithm~2.1 needs at most
\[
\left \lfloor \kappa_s ( f(x_0)- \flow)
\max\left( \epsilon_1^{-\frac{p+1}{p}}, \epsilon_2^{-\frac{p+1}{p-1}} \right)
\right \rfloor+1
\]
successful iterations (each involving one evaluation of $f$ and its $p$ first derivatives)
and at most
\[
\left \lfloor
\kappa_s ( f(x_0)- \flow)
\max\left( \epsilon_1^{-\frac{p+1}{p}}, \epsilon_2^{-\frac{p+1}{p-1}} \right)
\right \rfloor
                 \left(1+\frac{|\log\gamma_1|}{\log\gamma_2}\right)+
\frac{1}{\log\gamma_2}\log\left(\frac{\sigma_{\max}}{\sigma_0}\right)+1
\]
iterations in total to produce an iterate $x_\epsilon$ such that
$\|\nabla_x^1 f(x_\epsilon)\| \leq \epsilon_1$ and
$\lambdaleft\left(\nabla^2_x f(x_\epsilon)\right) \geq - \epsilon_2$,
where $\sigma_{\max}$ is given by \req{sigmaupper} and where
\[
\kappa_s \eqdef \frac{p+1}{\eta_1 \sigma_{\min}} \max\left(
\left(L + \theta +\sigma_{\max}\right)^{\frac{p+1}{p}},
\left((p-1)L + \theta +p\sigma_{\max}\right)^{\frac{p+1}{p-1}} \right).
\]
}

\proof{At each successful iteration $k$ before termination, 
either the first order or the second order approximate optimality condition must fail (at the next iteration), namely,
\eqn{fail-opt}{\chi_{f,1}(x_{k+1})>\epsilon_1\tim{or} \chi_{f,2}(x_{k+1})>\epsilon_2,}
and we also have the guaranteed decrease
\eqn{fdec}{
f(x_k)-f(x_{k+1})
\geq \eta_1 (T_p(x_k,0)-T_p(x_k,s_k))
\geq \bigfrac{\eta_1 \sigma_{\min}}{p+1} \;\|s_k\|^{p+1}
}
where we used \req{rhokdef2}, \req{Dphi} and \req{sigupdate2}. For any successful iteration for which the first condition in \req{fail-opt} holds, 
we deduce from \req{fdec}, \req{longs-g} and \req{sigmaupper} that
\eqn{eps1-decr}{
f(x_k)-f(x_{k+1}) \geq \kappa_1 \epsilon_1^{\frac{p+1}{p}}
\tim{where} \kappa_1 \eqdef
\bigfrac{\eta_1 \sigma_{\min}}{p+1}
\left(\frac{1}{L + \theta +\sigma_{\max}}\right)^{\frac{p+1}{p}}.}
Similarly, for any successful iteration for which the second condition in \req{fail-opt} holds, 
we deduce from \req{fdec}, \req{longs-H} and \req{sigmaupper} that
\eqn{eps2-decr}{
f(x_k)-f(x_{k+1}) \geq \kappa_2 \epsilon_2^{\frac{p+1}{p-1}}
\tim{where} \kappa_2 \eqdef \bigfrac{\eta_1 \sigma_{\min}}{p+1}
\left(\frac{1}{(p-1)L + \theta +p\sigma_{\max}}\right)^{\frac{p+1}{p-1}}.}
Thus on any successful iteration until termination we can guarantee the minimal of the two decreases in \req{eps1-decr} and \req{eps2-decr},
and hence,  since $\{f(x_k)\}$ decreases monotonically,
\[
f(x_0)-f(x_{k+1}) \geq \min[\kappa_1,\kappa_2] \min\left[\epsilon_1^{\frac{p+1}{p}}, \epsilon_2^{\frac{p+1}{p-1}}\right]\cdot|\mathcal{S}_k|.
\]
Using that $f$ is bounded below by $f_{\rm low}$, we conclude 
\[
| \calS_k | \leq \frac{f(x_0) - \flow}{\min[\kappa_1,\kappa_2]} \max\left[\epsilon_1^{-\frac{p+1}{p}}, \epsilon_2^{-\frac{p+1}{p-1}}\right]
\]
until termination,
from which the desired bound on the number of successful iterations follows.
Lemma~\ref{SvsU} is then invoked to compute the upper bound on the total
number of iterations.
}

Observe that we may modify the algorithm to seek only first-order points by
restricting \req{mterm} to $i=1$.
The corresponding complexity is then
\[
O\left(\epsilon_1^{-\frac{p+1}{p}}\right),
\]
which coincides with the bound in \cite{BirgGardMartSantToin17:mp}.
Moreover the same complexity result holds if, by chance,
$\lambdaleft\left(\nabla^2_x f(x_k)\right) \geq - \epsilon_2$ for all
iterations. By contrast, if $\epsilon_1$ is so large that
$\|\nabla_x^1 f(x_k)\| \leq \epsilon_1$ at every iteration, the complexity
is
\[
O\left(\epsilon_2^{-\frac{p+1}{p-1}}\right)
\]
to find a point with a sufficiently large leftmost eigenvalue.

\numsection{Final comments}
\label{sec-conclusion}

Our goal has been to devise an algorithm that can guaranteed to
find an approximate first- and second-order critical point in
fewer evaluations than the best known current champions. The new algorithm
we have designed finds such a point in at most
\[O \left(\max\left( \epsilon_1^{-\frac{p+1}{p}}, \epsilon_2^{-\frac{p+1}{p-1}}
\right) \right)\]
 function and derivative evaluations
under suitable differentiablity and Lipschitz continuity conditions. When $p=2$,
we recover the standard best bound \req{arc-eval}, while for $p=3$,
this improves to
$O \left(\max\left( \epsilon_1^{-4/3}, \epsilon_2^{-2} \right) \right)$
function and derivative evaluations, and approaches
$O \left(\max\left( \epsilon_1^{-1}, \epsilon_2^{-1} \right) \right)$
evaluations as $p$ increases to infinity. Of course, this comes at an
increased cost of requiring derivatives of order up to $p$, and
of needing to approximately solve a potentially harder step subproblem. Note though, that the conditions \req{descent2} and \req{mterm}
for model minimization are only local ones, and that the improved second-order approximate criticality result is achieved under the same
problem assumptions as the first order one (in \cite{BirgGardMartSantToin17:mp} and here).

In practice, the test \req{f-term} for termination in Step 1 of
Algorithm~2.1 would be arranged to check one of the pair
of required inequalities, and only to check the other if the
first holds (the order is immaterial). One could imagine a variant
of the algorithm in which failure of one (but not both) of \req{f-term}
might influence the requirement for the next step calculation/model minimization. Specifically,
if $\chi_{f,1}(x_k) > \epsilon_1$, one might simply require that
$\chi_{m,1}(x_k,s_k,\sigma_k) \leq \theta \|s_k\|^{p}$
rather than \req{mterm} as this alone would aim to improve first-order criticality. However, though this decoupling is possible
both in practice and in the analysis, it is not as straightforward as in the case of say, trust-region methods \cite{GratRoyoVice17}, as the 
lower bounds on the step in \req{longs-g} and \req{longs-H} depend on the objective's gradient and Hessian
value at the next trial point/iterate, not the current $x_k$.
Also, one might modify the ARp algorithm to check the optimality measures \req{f-term} at every trial point, not just successful ones.
This may allow earlier termination but possibly at an unsuccessful step and
at increased first- and second-derivatives evaluation cost.

While one might be tempted to try to provide bounds for an algorithm that
guarantees approximate third- and higher-order necessary optimality conditions,
we have not, as yet, been able to do so. 
The main sticking point has been that third-order necessary conditions
involve the behaviour of the third-order term of the Taylor series in
the nullspace of the Hessian (if it exists) 
and that this (typically proper) subspace of $\Re^n$ is highly sensitive.  Its use or the use of an approximating set is therefore open to miss-diagnosis.  
Higher-order criticality
becomes successively trickier; the critical spaces are then no longer
subspaces but cones \cite{CartGoulToin16a}.

Extending the approach here to the constrained case, even convex constraints, also seems challenging as the connection between model eigenvalues
and function eigenvalues in a set is no longer straightforward. Another aspect for future work is quantifying the cost of the subproblem solution
in a  similar vein to recent works \cite{Hazan, Duchi}, where there is particular interest due to large scale applications, in quantifying the number of derivative actions required per iteration
as derivatives cannot be stored/called explicitly. More generally, finding efficient ways to solve higher order polynomial models would bring ARp methods closer to practical use.


{\footnotesize

}


\newcommand{\resetcounters}{\setcounter{equation}{0} \setcounter{figure}{0}
 \setcounter{table}{0}}
\newcommand{\newappendixname}{A}
\newcommand{\newappendix}[1]{\renewcommand{\newappendixname}{{#1}}
 \section*{Appendix \newappendixname}
 \resetcounters
 \renewcommand{\theequation}{\newappendixname.\arabic{equation}}
 }

\renewcommand\appendix{%
  \setcounter{section}{0}%
  \setcounter{subsection}{0}%
  \renewcommand\thesection{\Alph{section}}}

\renewcommand{\thesection}{Appendix \Alph{section}}
\setcounter{section}{0}%
\setcounter{subsection}{0}%


\newappendix{A}
\renewcommand{\theequation}{A.\arabic{equation}}
\renewcommand{\thesection}{A.\arabic{section}}
\renewcommand{\thesubsection}{A.\arabic{subsection}}

\section{Proof of Lemma~\ref{taylor-bounds-lemma}}\label{A-Lemma-proof}

As in \cite{CartGoulToin16a}, consider the Taylor identity
\beqn{unitaylor}
\phi(1) - \tau_k(1) = \frac{1}{(k-1)!}\int_0^1 (1 - \xi)^{k-1}
                           [\phi^{(k)}(\xi ) - \phi^{(k)}(0)]  \, d\xi
\eeqn
involving a given univariate $C^k$ function $\phi(\alpha)$
and its $k$-th order Taylor
approximation
\[
\tau_k(\alpha) = \sum_{i=0}^k \phi^{(i)}(0) \frac{\alpha^i}{i!}
\]
expressed in terms of the value $\phi^{(0)} = \phi$ and
$i$th derivatives $\phi^{(i)}$, $i=1,\ldots,k$.
Then, picking $\phi(\alpha) = f(x+\alpha s)$ and $k = p$,
the identity
\beqn{intid}
\bigint_0^1 (1-\xi)^{k-1} \, d\xi = \frac{1}{k},
\eeqn
\req{tensor-Lip}, \req{Tnorm} and \req{unitaylor}
imply that, for all $x,s \in \Re^n$,
\[
f(x+s) \leq T_p(x,s) + \bigfrac{L}{p} \, \|s\|^{p+1}
\]
\cite[(2.8) with $L_{f,p} = (p-1)!L$]{CartGoulToin16a}
since $\tau_p(1) = T_p(x,s)$, which is the required \req{resf}.

Likewise, for an arbitrary unit vector $v$, selecting instead
$\phi(\alpha) = \nabla_x^1 f(x+\alpha s) [v]$ and $k=p-1$,
it follows from \req{unitaylor} that
\eqn{d1}{
\arr{rl}{
( \nabla_x^1  & \!\!\!\!\!\! f(x+ s) -  \nabla_s^1 T_p(x,s) ) [v] \\
=&\bigfrac{1}{(p-2)!}\bigint_0^1 (1 - \xi)^{p-2}
(\nabla_x^p f(x+ \xi s) - \nabla_x^p f(x))[s]^{p-1} [v] \, d\xi}
}
since $\tau_{p-1}(1) = \nabla_s^1 T_p(x,s)[v]$.
Thus, using the symmetry of the derivative tensors, picking  $v$ to
maximize the absolute value of the left-hand side of \req{d1} and using
\req{intid}, \req{Tnorm} and \req {tensor-Lip} successively, we
obtain that
\[
\begin{array}{ll}
\lefteqn{\!\!\!\!\!\| \nabla_x^1 f(x+ s) -  \nabla_s^1 T_p(x,s) \|_{[1]} }&
\\
& = \bigfrac{1}{(p-2)!}\left | \bigint_0^1 (1 - \xi)^{p-2}
(\nabla_x^p f(x+\xi s) - \nabla_x^p f(x)) [v]
\left[\frac{s}{\|s\|}\right]^{p-1}
\|s\|^{p-1}  d\xi \right |
\\*[2.5ex]
& \leq \bigfrac{1}{(p-2)!}\left[\bigint_0^1 (1 - \xi)^{p-2} d\xi\right]
\bigmax_{\xi \in [0,1]} \left | ( \nabla_x^p f(x+\xi s) - \nabla_x^p f(x) )
[v] \left[\frac{s}{\|s\|}\right]^{p-1} \right |
\|s\|^{p-1}
\\*[2.5ex]
& \leq \bigfrac{1}{(p-1)!} \bigmax_{\xi \in [0,1]}
\max_{\|w_1\|=\cdots=\|w_p\|=1} \left |
( \nabla_x^p f(x+\xi s) - \nabla_x^p f(x) ) [w_1,\ldots,w_p] \right |
\|s\|^{p-1}
\\*[2.5ex]
& = \bigfrac{1}{(p-1)!}
\bigmax_{\xi \in [0,1]}
\| \nabla_x^p f(x+\xi s) - \nabla_x^p f(x) \|_{[p]}
\|s\|^{p-1}
\\*[2.5ex]
& \leq L \|s\|^p
\end{array}
\]
which gives \req{resg}.

Finally, for arbitrary unit vectors $v_1$ and $v_2$, choosing
$\phi(\alpha) = \nabla_x^2 f(x+\alpha s) [v_1,v_2]$ and $k=p-2$,
the identity $\tau_{p-2}(1) = \nabla_s^2 T_p(x,s)[v_1,v_2]$ and
\req{unitaylor} together show that
\eqn{d2}{\arr{rl}{
 ( \nabla_x^2 & \!\!\!\!\!\! f(x+ s) -  \nabla_s^2 T_p(x,s) ) [v_1,v_2] \\
 =&\bigfrac{1}{(p-3)!}\bigint_0^1 (1 - \xi)^{p-3}
 (\nabla_x^p f(x+ \xi s) - \nabla_x^p f(x))[v_1,v_2] [s]^{p-2} \, d\xi.}}
As before, picking $v_1$ and $v_2$ to maximize the absolute value of the
left-hand side of \req{d2},
\[
\arr{rl}{
\| \nabla_x^2 & \!\!\!\!\!\! f(x+ s) -  \nabla_s^2 T_p(x,s) \|_{[2]} \\
= &\bigfrac{1}{(p-3)!} \left| \bigint_0^1 (1 - \xi)^{p-3}
(\nabla_x^p f(x+ \xi s) - \nabla_x^p f(x)) [v_1,v_2]
 \left[\frac{s}{\|s\|}\right]^{p-2} \|s\|^{p-2} \, d\xi \right| \\
\leq & \bigfrac{1}{(p-3)!} \left[ \bigint_0^1 (1 - \xi)^{p-3} d\xi \right]
\bigmax_{\xi \in [0,1]} \left|
(\nabla_x^p f(x+ \xi s) - \nabla_x^p f(x)) [v_1,v_2]
 \left[\frac{s}{\|s\|}\right]^{p-2} \|s\|^{p-2} \right| \\
\leq & \bigfrac{1}{(p-2)!}
\bigmax_{\xi \in [0,1]} \max_{\|w_1\|=\cdots=\|w_p\|=1} \left|
(\nabla_x^p f(x+ \xi s) - \nabla_x^p f(x)) [w_1,\ldots,w_p] \right
\|s\|^{p-2} \\
= & \bigfrac{1}{(p-2)!}
\bigmax_{\xi \in [0,1]}
\|\nabla_x^p f(x+ \xi s) - \nabla_x^p f(x)\|_{[p]} \|s\|^{p-2} \\
\leq & (p-1) L \|s\|^{p-1}
}\]
again using \req{tensor-Lip}, \req{Tnorm} and \req{intid}.
which provides \req{H-Lip}.


\section{Proof of Lemmas in Section \ref{sec-analysis}}

{\bf Proof of Lemma \ref{Dm_lemma}}
(See \cite[Lemma~2.1]{BirgGardMartSantToin17:mp})
Observe that, because of \req{descent2} and \req{mdef},
\[
0 < m(x_k,0,\sigma_k) - m(x_k,s_k,\sigma_k)
= T_p(x_k,0) - T_p(x_k,s_k) -\frac{\sigma_k}{p+1} \|s_k\|^{p+1}
\]
which implies the desired bound. Note that $s_k\neq 0$ as long as we can satisfy condition \req{descent2}, and so \req{Dphi} implies \req{rhokdef2} is well defined.\hfill$\Box$\\

{\bf Proof of Lemma \ref{sigmaupper_lemma}}
(See \cite[Lemma~2.2]{BirgGardMartSantToin17:mp})
Assume that
\beqn{siglarge}
\sigma_k \geq \frac{L(p+1)}{p \,(1-\eta_2)}.
\eeqn
Using \req{resf} and \req{Dphi}, we may then deduce that
\[
|\rho_k - 1|
\leq \frac{|f(x_k+s_k) - T_p(x_k,s_k)|}{|T_p(x_k,0)-T_p(x_k,s_k)|}
\leq \frac{L(p+1)}{p \,\sigma_k}
\leq 1-\eta_2
\]
and thus that $\rho_k \geq \eta_2$. Then iteration $k$ is very successful in
that $\rho_k \geq \eta_2$ and $\sigma_{k+1}\leq \sigma_k$.  As a consequence,
the mechanism of the algorithm ensures that \req{sigmaupper} holds.\hfill$\Box$\\

{\bf Proof of Lemma \ref{longs-g-lemma}}
(See \cite[Lemma~2.3]{BirgGardMartSantToin17:mp})
Using the triangle inequality,
\req{resg},
\req{gradm} and \req{mterm} for $i=1$, we obtain that
\[
\begin{array}{lcl}
\chi_{f,1}(x_k+s_k)
& \leq & \| \nabla_x^1 f(x_k+s_k) - \nabla_s^1 T_p(x_k,s_k)\|
        +\left\|\nabla_s^1 T_p(x_k,s_k) + \sigma_k\|s_k\|^p\bigfrac{s_k}{\|s_k\|}\right\| \\*[1.5ex]
 &       & + \sigma_k \|s_k\|^p\\*[1.5ex]
& =  & \| \nabla_x^1 f(x_k+s_k) - \nabla_s^1 T_p(x_k,s_k)\|_{[1]}+ \chi_{m,1}(x_k,s_k,\sigma_k) + \sigma_k \|s_k\|^p\\*[1.5ex]
 &\leq & L \|s_k\|^p + \chi_{m,1}(x_k,s_k,\sigma_k) + \sigma_k \|s_k\|^p\\*[1.5ex]
 & \leq & \left[L + \theta +\sigma_k \right] \|s_k\|^p
 \end{array}
\]
and \req{longs-g} follows.\hfill$\Box$\\

{\bf Proof of Lemma \ref{SvsU}.}
The regularization parameter update \req{sigupdate2}
gives that, for each $k$,
\[
\gamma_1\sigma_j \leq \max[\gamma_1\sigma_j,\sigma_{\min}]
 \leq \sigma_{j+1}, \ms j \in \calS_k,
\tim{ and }
\gamma_2\sigma_j \leq \sigma_{j+1}, \ms j \in \calU_k.
\]
Thus we deduce inductively that
\[
\sigma_0\gamma_1^{|\calS_k|}\gamma_2^{|\calU_k|}\leq \sigma_{k}.
\]
We therefore obtain, using \req{sigmax}, that
\[
|\calS_k|\log{\gamma_1}+|\calU_k|\log{\gamma_2}\leq
\log\left(\frac{\sigma_{\max}}{\sigma_0}\right),
\]
which then implies that
\[
|\calU_k|
\leq - |\calS_k|\frac{\log\gamma_1}{\log\gamma_2}
      + \frac{1}{\log\gamma_2}\log\left(\frac{\sigma_{\max}}{\sigma_0}\right),
\]
since $\gamma_2>1$. The desired result \req{unsucc-neg} then follows
from the equality $k +1 = |\calS_k| + |\calU_k|$
and the inequality $\gamma_1 < 1$ given by \req{eta-gamma2}.\hfill$\Box$


\begin{thebibliography}{10}

\bibitem{BirgGardMartSantToin17:mp}
E.G. Birgin, J.L. Gardenghi, J.M. Mart\'{i}nez, S.A. Santos, and Ph.~L. Toint.
\newblock Worst-case evaluation complexity for unconstrained nonlinear
  optimization using high-order regularized models.
\newblock {\em Mathematical Programming, Series~A}, 163(1):359--368, 2017.

\bibitem{BirginMartinez17}
E.G. Birgin and  J.M. Mart\'{i}nez.
\newblock The use of quadratic regularization with a cubic descent condition for unconstrained optimization.
\newblock {\em SIAM Journal on Optimization}, 27(2):1049--1074, 2017.

\bibitem{CartGoulToin11:mp}
C.~Cartis, N.~I.~M. Gould, and {\relax Ph}.~L. Toint.
\newblock Adaptive cubic regularisation methods for unconstrained optimization.
  {P}art {I}: motivation, convergence and numerical results.
\newblock {\em Mathematical Programming, Series~A}, 127(2):245--295, 2011.

\bibitem{CartGoulToin11b:mp}
C.~Cartis, N.~I.~M. Gould, and {\relax Ph}.~L. Toint.
\newblock Adaptive cubic regularisation methods for unconstrained optimization.
  {P}art {II}: worst-case function and derivative-evaluation complexity.
\newblock {\em Mathematical Programming, Series~A}, 130(2):295--319, 2011.

\bibitem{CartGoulToin12d}
C.~Cartis, N.~I.~M. Gould, and {\relax Ph}.~L. Toint.
\newblock Complexity bounds for second-order optimality in unconstrained
  optimization.
\newblock {\em Journal of Complexity}, 28:93--108, 2012.

\bibitem{CartGoulToin12:opt}
C.~Cartis, N.~I.~M. Gould, and {\relax Ph}.~L. Toint.
\newblock How much patience do you have? a worst-case perspective on smooth
  nonconvex optimization.
\newblock {\em Optima}, 88:1--10, 2012.

\bibitem{CartGoulToin16a}
C.~Cartis, N.~I.~M. Gould, and {\relax Ph}.~L. Toint.
\newblock Second-order optimality and beyond: characterization and evaluation
  complexity in nonconvex convexly-constrained optimization.
\newblock Preprint RAL-P-2016-008, Rutherford Appleton Laboratory, Chilton,
  Oxfordshire, England, 2016.

\bibitem{ConnGoulToin00}
A.~R. Conn, N.~I.~M. Gould, and {\relax Ph}.~L. Toint.
\newblock {\em Trust-{R}egion {M}ethods}.
\newblock SIAM, Philadelphia, 2000.

\bibitem{CurtRobiSame17:mp}
F.~E. Curtis, D.~P. Robinson, and M.~Samadi.
\newblock A trust region algorithm with a worst-case iteration complexity of
  {$O(\epsilon^{-3/2})$} for nonconvex optimization.
\newblock {\em Mathematical Programming, Series~A}, 162(1):1--32, 2017.

\bibitem{GrapYuanYuan16:jota}
G.N. Grapiglia, J.~Yuan, and Y.~Yuan.
\newblock Nonlinear stepsize control algorithms: Complexity bounds for first
  and second-order optimality.
\newblock {\em Journal of Optimization Theory and Applications},
  17(3):980--997, 2016.

\bibitem{Duchi}
Y. Carmon and J. C. Duchi.
\newblock Gradient descent efficiently finds the cubic-regularized nonconvex Newton step.
\newblock Stanford University, Technical Report. ArXiv:1612.00547, 2016.

\bibitem{GratRoyoVice17}
S.~Gratton, C.~W. Royer, and L.~N. Vicente.
\newblock A decoupled first/second-order steps technique for nonconvex
  nonlinear unconstrained optimization with improved complexity bounds.
\newblock Technical Report 17-21, Department of Mathematics, University of
  Coimbra, 2017.

\bibitem{Hazan}
N. Agarwal, Z. Allen-Zhu, B. Bullins, E. Hazan and T. Ma.
\newblock Finding approximate local minima faster than gradient descent.
\newblock Princeton University, Technical Report, ArXiv: 1611.01146, 2016.


\bibitem{MartRayd16}
J.~M. Mart\'{i}nez and M.~Raydan.
\newblock Cubic-regularization counterpart of a variable-norm trust-region
  method for unconstrained minimization.
\newblock {\em Journal of Global Optimization}, 2016.
\newblock DOI: 10.1007/s10898-016-0475-8.

\bibitem{Nest04}
Y.~Nesterov.
\newblock {\em Introductory lectures on convex optimization}.
\newblock Kluwer Academic Publishers, Dordrecht, The Netherlands, 2004.

\bibitem{NestPoly06}
{Yu}. Nesterov and B.~T. Polyak.
\newblock Cubic regularization of {N}ewton method and its global performance.
\newblock {\em Mathematical Programming, Series~A}, 108(1):177--205, 2006.


\end{thebibliography}
\end{document}